\documentclass{amsart}

\usepackage{amssymb}

\newtheorem{theorem}{Theorem}[section]
\theoremstyle{definition}
\newtheorem{prop}[theorem]{Proposition}
\newtheorem{lemma}[theorem]{Lemma}
\newtheorem{cor}[theorem]{Corollary}
\newtheorem{definition}[theorem]{Definition}

\numberwithin{equation}{section}

\begin{document}

\newcommand{\cc}{2^{\aleph_0}}
\newcommand{\Bb}{{\mathbb B}}
\newcommand{\Cb}{{\mathbb C}}
\newcommand{\Ac}{{\mathcal A}}
\newcommand{\Bc}{{\mathcal B}}
\newcommand{\Cc}{{\mathcal C}}
\newcommand{\Dc}{{\mathcal D}}
\newcommand{\Hc}{{\mathcal H}}
\newcommand{\Mca}{{\mathcal M}}
\newcommand{\Kc}{{\mathcal K}}
\newcommand{\Pc}{{\mathcal P}}
\newcommand{\Xc}{{\mathcal X}}
\newcommand{\fin}{{\rm fin}}
\newcommand{\omf}{P(\omega)/\fin}
\newcommand{\subs}{\subseteq}

\title{$P(\omega)/{\rm fin}$ and projections in the Calkin algebra}

\author{Eric Wofsey}
\address{Department of Mathematics, Washington University in Saint Louis, Saint Louis, MO 63130}
\email{erwofsey@artsci.wustl.edu}


\date{August 29, 2006, and, in revised form, December 28, 2006.}

\commby{Julia Knight}

\begin{abstract}
We investigate the set-theoretic properties of the lattice of projections in the Calkin algebra of a separable infinite-dimensional Hilbert space in relation to those of the Boolean algebra $\omf$, which is isomorphic to the sublattice of diagonal projections.  In particular, we prove some basic consistency results about the possible cofinalities of well-ordered sequences of projections and the possible cardinalities of sets of mutually orthogonal projections that are analogous to well-known results about $\omf$.
\end{abstract}

\maketitle

We consider the Hilbert space $\Hc=l^2=l^2(\omega)$ with standard basis $\{e_n\}$ and we denote the algebra of bounded operators on $\Hc$ by $\Bc$.  Let $\Kc\subset\Bc$ be the ideal of compact operators; that is, the norm-closure of the ideal of finite-rank operators (for background on the elementary properties of compact operators, see \cite{Con}).  Let $\Cc$ be the quotient $\Bc/\Kc$ (the Calkin algebra), $\pi:\Bc\rightarrow\Cc$ be the quotient map, and $\Pc$ be the lattice of projections in $\Cc$.  For $A\subs\omega$, we let $P_A$ be the projection onto $l^2(A)\subs l^2(\omega)$; the map $A\mapsto P_A$ embeds the Boolean algebra $P(\omega)$ into the lattice of projections in $\Bc$.  Let $\fin$ be the ideal of finite subsets of $\omega$; we will not distinguish between a subset $A\subs\omega$ and its coset under the quotient map $P(\omega)\rightarrow P(\omega)/\fin$.  Then $A \mapsto \pi(P_A)$ is, in fact, a well-defined embedding of $P(\omega)/\fin$ into $\Pc$, the \emph{diagonal embedding}.  The diagonal embedding preserves the lattice operations and takes Boolean complements to orthogonal complements.

In \cite{Had}, D.\ Hadwin showed that under the Continuum Hypothesis, all maximal chains in $\Pc$ are order-isomorphic, and conjectured that this condition was equivalent to CH.  The methods used in that proof were essentially the same as those that can be used to prove similar and related statement about $\omf$ under CH.  In this paper, we will expand on this parallel, showing how projections in $\Bc$ can be considered to be analogous to subsets of $\omega$, and how the quotient $\Cc=\Bc/\Kc$ is then analogous to $\omf$.  We do not settle Hadwin's conjecture, but do show (Corollary \ref{HaCo}) that it is consistent for non-isomorphic maximal chains in $\Pc$ to exist.

Two of the basic problems concerning $\omf$ are:
\begin{enumerate}
\item What kind of maximal well-ordered sequences (``limits'') are there in the partial ordering of $\omf$?
\item What cardinalities can an (infinite) maximal set of disjoint elements of $\omf$ (a \emph{maximal almost disjoint family}, or maximal adf) have?
\end{enumerate}
It is easy to show that any limit or maximal adf must be uncountable.  In the presence of the continuum hypothesis (in fact, Martin's Axiom suffices), both of these problems are thus trivial; in general, neither problem can be decided in ZFC.  It is, however, well-known that there exists an adf of cardinality $\cc$, which by Zorn's Lemma extends to a maximal adf of cardinality $\cc$.

If we ``apply the diagonal embedding'' to these two problems, we obtain:
\begin{enumerate}
\item What kind of limits are there in the partial ordering of $\Pc$?
\item What cardinalities can an (infinite) maximal set of orthogonal elements of $\Pc$ (a \emph{maximal almost orthogonal family}, or maximal aof) have?
\end{enumerate}
Just as in the case of $\omf$, these problems are trivial under the continuum hypothesis (see \cite{Had} for the first problem); in Section 3 we show that they are also easily solved with Martin's Axiom.  Every well-ordered sequence in $\omf$ clearly maps to a well-ordered sequence in $\Pc$ under the diagonal embedding, and similarly any adf maps to an aof.  Since there exists a (maximal) adf of cardinality $\cc$, the same is true of aofs.  A natural question is whether limits (i.e., maximal sequences) and maximal adfs remain maximal under the diagonal embedding.  In Section 1 we show that certain generic sequences and maximal adfs do remain maximal, while in Section 2 we show that in general, limits and maximal adfs do not remain maximal.  We assume in Sections 1 and 3 that the reader is familiar with the basic language and methods of forcing; see \cite{Kun} for a good introduction to the subject.

The questions discussed in this paper are only a few of the set-theoretic questions concerning the Calkin algebra and $\Pc$.  More generally, one could ask whether the diagonal embedding can be used to find analogs of the classical cardinal invariants related to $\omf$ in the setting of $\Pc$, and whether these ``quantized'' cardinal invariants have the same values of the classical ones.  For a discussion of other set-theoretic problems about the Calkin algebra that are of interest from a more analytic perspective, see \cite{Wea}.

I would like to thank Nik Weaver for his supervision, guidance, advice, and contributions to this paper, and for teaching me all of the background I needed.  I also thank Ilijas Farah for some helpful suggestions.

\section{Forcing limits in $\Pc$ and maximal aofs}
First, we prove a useful characterization of the partial ordering on $\Pc$ (noting that $\pi(p)\leq\pi(q)$ iff $p(1-q)$ is compact).

\begin{definition}
Let $p$ be a projection in $\Bc$, let $\{a_n\}$ be an orthonormal subset of $\Hc$, and let $\epsilon>0$.  Then an \emph{$\epsilon$-block} for $p$ (with respect to $\{a_n\}$) is a pair $(S,v)$ such that $S$ is a finite subset of $\omega$, $v\in{\rm ran}(p)$ is a unit vector, and $\|P_S(v)\|>\epsilon$ (where $P_S$ is the projection onto $\overline{{\rm span}}\{a_n:n\in S\}$).  We say $(S,v)$ is an $\epsilon$-block \emph{above $N$} if for each $m\in S$, $m>N$.
\end{definition} 

\begin{lemma}\label{blo}
Let $p$ and $q$ be projections in $\Bc$, and let $\{a_n\}$ be an orthonormal basis for ${\rm ran}(q)$.  Then $pq$ is compact iff $\forall\epsilon>0$, $\exists N\in\omega$ such that there are no $\epsilon$-blocks for $p$ with respect to $\{a_n\}$ above $N$.  In particular, $p$ is compact iff this holds with $\{a_n\}=\{e_n\}$, the standard basis, and for $A\subs\omega$ and $q=P_A$, $pq$ is compact iff $\forall\epsilon>0$, $\exists N\in\omega$ such that there are no $\epsilon$-blocks $(S,v)$ for $p$ with respect to $\{e_n\}$ above $N$ such that $S\subs A$.
\end{lemma}
\begin{proof}
($\Leftarrow$):  Suppose the second condition holds; fix $\epsilon>0$ and $N$ such that there are no $\epsilon$-blocks above $N$.  Let $q_0$ be the projection onto $\overline{{\rm span}}\{a_n:n\leq N\}$ and $q'=q-q_0$.  Since $q_0$ is compact,
  \[\|\pi(pq)\|=\|\pi(pq')\|=\|\pi((pq')^*)\|=\|\pi(q'p)\|\leq\|q'p\|\leq\epsilon,\]
with the final inequality holding by our choice of $N$.  Since $\epsilon$ is arbitrary, $\pi(pq)=0$.
\\\\
($\Rightarrow$):  Let $q_n$ be the projection onto $\overline{{\rm span}}\{a_k:k>n\}$, and suppose the second condition fails but $pq$ is compact.  Then choose $\epsilon>0$ and a sequence of unit vectors $v_n$ in ${\rm ran}(p)$ such that $\|q_n(v_n)\|>\epsilon$ for all $n$.  Since $pq$ is compact, $(pq)^*=qp$ is too, so the image of the unit ball under $qp$ is precompact.  But each $q(v_n)=qp(v_n)$ is in that image, so some subsequence $(q(v_{n_k}))$ converges.  For any $k$, let $N(k)>k$ be such that $\|q_{n_N(k)}(v_{n_k})\|<\epsilon/2$; then
  \[\|q(v_{n_k})-q(v_{n_N(k)})\|\geq\|q_{n_N(k)}(v_{n_k})-q_{n_N(k)}(v_{n_N(k)})\|\geq \epsilon/2.\]
Since this holds for all $k$, $(q(v_{n_k}))$ cannot be a Cauchy sequence, a contradiction. 
\end{proof}
\noindent All $\epsilon$-blocks that we consider will be with respect to $\{e_n\}$ unless stated otherwise.

Now let us consider forcing limits in $\omf$ and $\Pc$.

\begin{definition}
Let $L$ be a partially ordered set and $\kappa$ be a regular cardinal.  A \emph{$\kappa$-limit} in $L$ is a pair $((a_\alpha)_{\alpha<\kappa},b)$ with $a_\alpha,b\in L$ such that $(a_\alpha)$ is strictly increasing, $a_\alpha<b$ for all $\alpha$, and there is no $c\in L$ such that $a_\alpha<c<b$ for all $\alpha$.  If such a $c$ does exist, we say $c$ \emph{interpolates} $((a_\alpha),b)$.\label{gap}
\end{definition}

A simple forcing argument (originally from \cite{Hec}) shows that for any regular $\kappa\leq\cc$, it is consistent for there to be $\lambda$-limits in $\omf$ for all regular uncountable $\lambda\leq\kappa$.  We will show that the generic limit added by that notion of forcing remains maximal under the diagonal embedding.

Let $\alpha$ be an ordinal.  We define
  \[T_\alpha=\{f:F\times n\rightarrow 2: F\textnormal{ is a finite subset of $\alpha$ and $n\in\omega\}$}.\]
For $f:F\times n\rightarrow 2$ and $g:F'\times n'\rightarrow 2$ in $T_\alpha$, we say $g\leq f$ if $g\supseteq f$ and for all $\beta,\gamma\in F$ with $\beta<\gamma$ and for all $k$ such that $n\leq k<n'$, $g(\beta,k)\leq g(\gamma,k)$.  For $\alpha<\beta$, $T_\alpha$ is clearly completely embedded in $T_\beta$.  For $f:F\times n\rightarrow 2$ and $g:F'\times n\rightarrow 2$ in $T_\alpha$, if $f\cup g$ is a function then it is a common extension of $f$ and $g$ in $T_\alpha$.  For any uncountable $\{f_\beta:F_\beta\times n_\beta\rightarrow 2\}\subs T_\alpha$, $n_\beta$ is constant on an uncountable subset, so by the $\Delta$-system lemma $T_\alpha$ is ccc.

\begin{theorem}
Let $\kappa$ be a regular uncountable cardinal.  Then $T_\kappa\Vdash$``for every regular uncountable cardinal $\lambda\leq\kappa$ there is a $\lambda$-limit in $\Pc$.''  Also, if $\kappa^{\aleph_0}=\lambda$ (in the base model), $T_\kappa\Vdash$``$\cc\leq\lambda$.''\label{gf}
\end{theorem}
\begin{proof}
First, note that the bound on $\cc$ follows from the fact that $T_\kappa$ is ccc and has cardinality $\kappa$.  Let $\Mca$ be a transitive model of a sufficiently large fragment of ZFC with $\kappa\in\Mca$, and let $\Mca_0\subs \Mca$ be a countable elementary submodel containing $\kappa$.  Denote the transitive collapse of $\Mca_0$ by $M$; we will work in $M$ and not distinguish between elements of $\Mca_0$ and the corresponding elements of $M$.

Let $G\subset T_\kappa$ be an $M$-generic filter.  For $\alpha\leq\kappa$ and $n<\omega$, write $g_\alpha(n)=(\bigcup G)(\alpha,n)$, $G_\alpha=G\cap T_\alpha$, and $M_\alpha=M[G_\alpha]$.  Define $A_\alpha=g_\alpha^{-1}(\{1\})$; it is easy to see that the sets $A_\alpha$ form an increasing sequence in $\omf$.  We claim that, in fact, $((\pi(P_{A_\alpha})),1)$ is a $\kappa$-limit in $\Pc$.

Let $p\in M_\kappa$ be a projection in $\Bc$, and suppose $\pi(p)$ interpolates $((\pi(P_{A_\alpha})),1)$.  Then in particular, $1-p$ is not compact, so there exist $\epsilon>0$ and sequences $(S_n)$ and $(v_n)$ such that for each $n$, $(S_n,v_n)$ is an $\epsilon$-block for $1-p$ above $n$.  By considering $p$ as a matrix of countably many reals and each $v_n$ as a sequence of reals, it is clear that there is a definable bijection between reals and such triples $(p,(S_n),(v_n))$.  Since $T_\kappa$ is ccc, by taking nice names we see that every real in $T_\kappa$ has a $T_\alpha$-name for some $\alpha<\kappa$.  In particular, we choose a nice name for the real associated to $(p,(S_n),(v_n))$ to obtain an $\alpha<\kappa$ such that $(p,(S_n),(v_n))\in M_\alpha$.  Let $\tilde{p}$, $(\tilde{S}_n)$, and $(\tilde{v}_n)$ be $T_\alpha$-names for $p$, $(S_n)$, and $(v_n)$.

Let $q\in G_\alpha$ force ``$\tilde{p}$ is a projection in $\Bc$, and $(\tilde{S}_n,\tilde{v}_n)$ is an $\epsilon$-block for $1-\tilde{p}$ above $n$ for each $n$''.  Define $D_N\subs T_\kappa$ as follows: $f:F\times n\rightarrow2$ is in $D_N$ iff there is some $S\subs n$ and some $T_\kappa$-name $\tilde{v}$ such that $f(\alpha,m)=1$ for all $m\in S$ and $f\Vdash$``$(S,\tilde{v})$ is an $\epsilon$-block for $1-\tilde{p}$ above $N$''.  We claim that for each $N$, $D_N$ is dense below $q$.  Indeed, let $f:F\times n\rightarrow 2$ extend $q$; it is no loss of generality to suppose $n<N$ and $\alpha\in F$.  % because an $\epsilon$-block above $\max(n,N)$ is also an $\epsilon$-block above $N$.
Since $q\Vdash$``$\tilde{S}_N$ is a finite subset of $\omega$'', we can extend $f$ to $f':F'\times n'\rightarrow 2$ such that $f'$ determines the value of $\tilde{S}_N$ and forces $\tilde{S}_N\subs n'$.  We may also assume that $f'(\beta,m)=1$ for every $m$ such that $f'\Vdash$``$m\in\tilde{S}_N$'' and for every $\beta\geq\alpha$ in $F$: modifying $f'$ such that this is true will not change the fact that $f'$ extends $f$ because $N>n$, and it will not change the value of $\tilde{S}_N$ since $\tilde{S}_N$ is a $T_\alpha$-name, and $f'$'s projection onto $T_\alpha$ is unchanged.  It is thus easy to see that $f'$ is an extension of $f$ in $D_N$, as desired.

Since $q\in G$ and each $D_N$ is dense below $q$, $G$ meets each $D_N$.  Hence for any $N$, there exists $v$ and $S\subs A_\alpha$ such that $(S,v)$ is an $\epsilon$-block above $N$ for $1-p$.  By Lemma \ref{blo}, $(1-p)P_{A_\alpha}$ is thus not compact, contradicting the assumption that $\pi(p)$ interpolates $((\pi(P_{A_\alpha})),1)$.

Hence $((\pi(P_{A_\alpha})),1)$ is a $\kappa$-limit in $\Pc$.  A similar argument shows that for any regular uncountable $\lambda<\kappa$, $((\pi(P_{A_\alpha}))_{\alpha<\lambda},\pi(P_{A_\lambda}))$ is a $\lambda$-limit.  Hence $T_\kappa\Vdash$``for every regular uncountable cardinal $\lambda\leq\kappa$ there is a $\lambda$-limit in $\Pc$.''
\end{proof}

\begin{cor}
``There exists a $\kappa$-limit in $\Pc$ for some $\kappa<2^{\aleph_0}$'' is consistent with ZFC.
\end{cor}

As another corollary, we obtain a result related to a conjecture in \cite{Had}.

\begin{cor} \label{HaCo}
``There exist non-isomorphic maximal chains in $\Pc$'' is consistent with ZFC.
\end{cor}

\begin{proof}
Start with any model of ZFC and force with $T_{\omega_2}$.  We then obtain $\omega_1$- and $\omega_2$-limits $((A_\alpha),1)$ and $((B_\beta),1)$ in $\Pc$ (if $((C_\alpha),D)$ is an $\omega_1$-limit, so is $((C_\alpha+(1-D)),1)=(A_\alpha,1)$).  Extend $\{A_\alpha\}$ and $\{B_\beta\}$ to maximal chains $L$ and $L'$.  Then any isomorphism from $L$ to $L'$ must take $1$ to $1$ and hence take $(A_\alpha)$ to a sequence that is cofinal with $(B_\beta)$, which is impossible.  Hence $L$ and $L'$ are not isomorphic.
\end{proof}

We now consider the analogous problem maximal aofs (again, the case of maximal adfs is from \cite{Hec}).  For $\Xc$ a set of disjoint sets, we define
  \[U_\Xc=\{f:F\times n\rightarrow 2: F\textnormal{ is a finite subset of $\bigcup\Xc$ and $n\in\omega\}$}.\]
For $f:F\times n\rightarrow 2$ and $g:F'\times n'\rightarrow 2$ in $U_\Xc$, we say $g\leq f$ if $g\supseteq f$ and for any $X\in\Xc$, for any distinct $x,y\in X\cap F$, and for all $k$ such that $n\leq k<n'$, $g(x,k)$ and $g(y,k)$ are not both $1$.  For $Y\subs\bigcup\Xc$, $U^Y_\Xc=U_{\{Y\cap X\}_{X\in\Xc}}$ is clearly completely embedded in $U_\Xc$.  By the same argument as for $T_\alpha$, $U_\Xc$ is ccc for any $\Xc$.

\begin{theorem}
Let $\Xc$ be set of disjoint uncountable sets.  Then $U_\Xc\Vdash$``for each $X\in\Xc$, there is a maximal aof of cardinality $|X|$.''  Also, if $|\bigcup\Xc|^{\aleph_0}=\lambda$ (in the base model), $U_\Xc\Vdash$``$\cc\leq\lambda$.''
\end{theorem}
\begin{proof}
The proof is almost identical to that of Theorem \ref{gf}.  First, the bound on $\cc$ follows from the fact that $U_\Xc$ is ccc.  We work in a countable transitive model $M$ as in Theorem \ref{gf}, and use the same notation (writing $G_Y=G\cap U^Y_\Xc$ for $Y\subs\bigcup\Xc$, and $M_Y=M[G_Y]$); we wish to show that for each $X\in\Xc$, $\{\pi(P_{A_x})\}_{x\in X}$ is maximal as an aof.  Fix $X\in\Xc$ and suppose $p\in M_{\bigcup\Xc}$ is a non-compact projection in $\Bc$ that is almost orthogonal to every $P_{A_x}$.  Then since $p$ is not compact, there exists $\epsilon>0$ and sequences $(S_n)$ and $(v_n)$ such that for each $n$, $(S_n,v_n)$ is an $\epsilon$-block for $p$ above $n$.  As in Theorem \ref{gf}, there is a countable set $Y\subset\bigcup\Xc$ such that $p$, $(S_n)$, and $(v_n)$ have $U^Y_\Xc$-names $\tilde{p}$, $(\tilde{S}_n)$, and $(\tilde{v}_n)$.  We fix $x_0\in X-Y$.

Let $q\in G_Y$ force ``$\tilde{p}$ is a projection in $\Bc$, and $(\tilde{S}_n,\tilde{v}_n)$ is an $\epsilon$-block for $\tilde{p}$ above $n$ for each $n$''.  Define $D_N\subs U_\Xc$ as follows: $f:F\times n\rightarrow2$ is in $D_N$ iff there is some $S\subs n$ and some $U_\Xc$-name $\tilde{v}$ such that $f(x_0,m)=1$ for all $m\in S$ and $f\Vdash$``$(S,\tilde{v})$ is an $\frac{\epsilon}{\sqrt{2}}$-block for $\tilde{p}$ above $N$''.  We claim that for each $N$, $D_N$ is dense below $q$.  Indeed, let $f:F\times n\rightarrow 2$ extend $q$; it is no loss of generality to suppose $n<N$ and $x_0\in F$.  Let $B=\bigcup_{x\in F\cap X-\{x_0\}} A_x$ and let $\tilde{B}$ be its canonical name.  Since $p$ is almost orthogonal to each $A_x$ (i.e. $pP_{A_x}$ is compact), it is almost orthogonal to $B$.  Therefore we can extend $f$ to $f'$ such that for some $M$, $f'\Vdash$``there are no $\frac{\epsilon}{\sqrt{2}}$-blocks $(S,v)$ for $\tilde{p}$ above $M$ such that $S\subs \tilde{B}$''; we lose no generality by assuming $M<N$.  Since $q\Vdash$``$\tilde{S}_N$ is a finite subset of $\omega$'', we can extend $f'$ to $f'':F''\times n''\rightarrow 2$ such that $f''$ determines the value of $\tilde{S}_N$ and forces $\tilde{S}_N\subs n''$.  Let
  \[\textnormal{$S=\{m:f''\Vdash$``$m\in\tilde{S}_N$ and $m\not\in \tilde{B}$''$\}$}\]
and
  \[\textnormal{$S'=\{m:f''\Vdash$``$m\in\tilde{S}_N$ and $m\in \tilde{B}$''$\}$.}\]
Since $f''$ extends $f'$, $f''$ forces that $(S',\tilde{v}_N)$ is not an $\frac{\epsilon}{\sqrt{2}}$-block for $\tilde{p}$.  Thus since $f''$ forces $(S\cup S',\tilde{v}_N)$ to be an $\epsilon$-block, it also forces $(S,\tilde{v}_N)$ to be an $\frac{\epsilon}{\sqrt{2}}$-block.

Now let $h$ be $f''$ modified such that $h(x_0,m)=1$ for every $m\in S$.  Then $h$ still extends $f$ since $N>n$ and by the definition of $S$, $f''(x,m)=0$ for $m\in S$ and $x\in F\cap X-\{x_0\}$.  Also, $h$ still forces $(S,\tilde{v}_N)$ to be an $\frac{\epsilon}{\sqrt{2}}$-block for $\tilde{p}$ since $\tilde{p}$ and $\tilde{v}_N$ are $U^Y_\Xc$-names and $f''$ and $h$ have the same projection onto $U^Y_\Xc$.  Thus $h$ is an extension of $f$ in $D_N$.

By genericity, then, for any $N$ there exist $\frac{\epsilon}{\sqrt{2}}$-blocks $(S,v)$ for $p$ above $N$ with $S\subs A_{x_0}$.  But then $p$ is not almost orthogonal to $P_{A_{x_0}}$, a contradiction.  Hence $\{\pi(P_{A_x})\}_{x\in X}$ is a maximal aof of cardinality $|X|$ for each $X\in\Xc$.
\end{proof}

\begin{cor}
``There exists a maximal aof of cardinality $<2^{\aleph_0}$'' is consistent with ZFC.
\end{cor}

\section{Counterexamples}
In this section we show a way of constructing limits in $\omf$ or maximal adfs that do not remain maximal under the diagonal embedding.  The key tool is the following:

\begin{definition}
A \emph{thinness criterion} is a partition of $\omega$ into finite sets $B_n$ such that $|B_n|\rightarrow\infty$ as $n\rightarrow\infty$.  A subset $A$ of $\omega$ is \emph{thin} (with respect to $(B_n)$) if $|A\cap B_n|/|B_n|\rightarrow 0$ as $n\rightarrow\infty$.  For $B$ a finite subset of $\omega$, let $1_B$ be the characteristic function of $B$ in $\omega$, considered as an $l^2$ sequence.  Then the \emph{dominating projection} $D_\Bb$ associated with a thinness criterion $\Bb=(B_n)$ is the projection onto $({\rm span}_n\{1_{B_n}\})^\perp$.
\end{definition}

\begin{prop}\label{thin}
Let $A\subs\omega$ and $\Bb$ be a thinness criterion.  Then $A$ is thin with respect to $\Bb$ iff $\pi(P_A)\leq\pi(D_\Bb)$.
\end{prop}
\begin{proof}
($\Rightarrow$): Suppose $A$ is thin; fix $\epsilon>0$ and choose $N$ such that $|A\cap B_n|/|B_n|<\epsilon$ for $n>N$.  Define $A_N=A-\bigcup_{n\leq N} B_n$.  Let $v\in{\rm ran}(P_{A_N})$ be a unit vector and write $v=\sum a_n v_n$ for $a_n\in\Cb$ and $v_n\in{\rm ran}(P_{A_N\cap B_n})$ unit vectors.  Then by our choice of $N$, 
  \[\|(1-D_\Bb)(v_n)\|^2=\left|\left\langle \frac{1_{B_n}}{\sqrt{|B_n|}},v_n\right\rangle\right|^2<\epsilon.\]
Thus
  \[\|(1-D_\Bb)(v)\|^2=\sum |a_n|^2 \|(1-D_\Bb)(v_n)\|^2<\epsilon\sum |a_n|^2=\epsilon.\]
Since $v$ was arbitrary, $\|\pi(1-D_\Bb)\pi(P_A)\|=\|\pi(1-D_\Bb)\pi(P_{A_N})\|\leq\|(1-D_\Bb)P_{A_N}\|<\epsilon$.  Since $\epsilon$ was arbitrary, $\pi(P_A)\leq\pi(D_\Bb)$.
\\\\
($\Leftarrow$):  Suppose $\pi(P_A)\leq\pi(D_\Bb)$ and fix $\epsilon>0$.  Then there is some $N$ beyond which there are no $\epsilon$-blocks for $P_A$ with respect to the basis $\{\frac{1_{B_n}}{\sqrt{|B_n|}}\}$.  Let $A_n=A\cap B_n$; then 
  \[\left\langle \frac{1_{B_n}}{\sqrt{|B_n|}}, \frac{1_{A_n}}{\sqrt{|A_n|}}\right\rangle=\sqrt{\frac{|A_n|}{|B_n|}}.\]
For $n>N$, $(\{n\},\frac{1_{A_n}}{\sqrt{|A_n|}})$ is not an $\epsilon$-block, so $|A_n|/|B_n|\leq\epsilon^2$.  Thus $|A_n|/|B_n|$ converges to $0$, so $A$ is thin.
\end{proof}

\begin{prop}
There exists a maximal adf that does not remain maximal under the diagonal embedding.
\end{prop}
\begin{proof}
Let $\Bb$ be a thinness criterion and let $f:[\omega]^\omega\rightarrow[\omega]^\omega$ take each infinite subset $S$ of $\omega$ to an infinite thin subset of $S$.  Now let $\Dc$ be a maximal almost disjoint subfamily of the range of $f$ (i.e. a maximal subset of $\{T: \exists S$ such that $f(S)=T\}$ that is an adf).  Then $\Dc$ is a maximal adf: if $S\subs\omega$ is infinite and almost disjoint from every element of $\Dc$, then $f(S)$ has the same property.  But then $\Dc\cup\{f(S)\}$ is an adf, contradicting the maximality of $\Dc$.  However, by Proposition \ref{thin}, $\pi(1-D_\Bb)$ is orthogonal to the image of every element of $\Dc$ under the diagonal embedding.  Hence $\Dc$ does not remain maximal under the diagonal embedding.
\end{proof}

\begin{prop}
CH implies that there is a limit in $\omf$ that does not remain a limit under the diagonal embedding.
\end{prop}
\begin{proof}
Let $\Bb=(B_n)$ be a thinness criterion and enumerate the coinfinite subsets of $\omega$ as $(S_\alpha)_{\alpha<\omega_1}$.  Define an increasing (mod ${\rm fin}$) sequence $(A_\alpha)_{\alpha<\omega_1}$ of thin sets by induction: let $A_0=\emptyset$.  If $\alpha$ is a limit ordinal, let $A_\alpha$ be a thin set such that $A_\beta\leq A_\alpha$ (mod $\fin$) $\forall\beta<\alpha$ (a simple diagonalization argument shows that this is possible since ${\rm cf}(\alpha)=\omega$).  If $\alpha=\beta+1$, let $T=\{k: k$ is the least element of $B_n\cap(\omega-S_\beta)$ for some $n\}$; since $S_\beta$ is coinfinite, $T$ is infinite.  Define $S_\alpha=S_\beta\cup T$; then $S_\alpha$ is still thin since $T$ contains at most one point from each $B_n$.

Now $((A_\alpha),1)$ is a limit in $\omf$: if $S=S_\alpha\subs\omega$ is coinfinite, then $A_{\alpha+1}\not\leq S$ (mod $\fin$).  But in $\Pc$, $\pi(D_\Bb)\geq\pi(P_{A_\alpha})$ for each $\alpha$ by Proposition \ref{thin}, so the image under the diagonal embedding is interpolated by $\pi(D_\Bb)$.
\end{proof}

\section{MA and $\Pc$}
It is well known that Martin's Axiom implies that all (infinite) maximal adfs have cardinality $\cc$ and that $\kappa$-limits exist in $\omf$ only for $\kappa=\cc$ (see \cite{Kun} for a proof for maximal adfs; the proof for limits is similar).  We now prove the analogous results for $\Pc$.

\begin{lemma}\label{malem}
Let $p$ and $q$ be projections in $\Bc$ and let $\{a_n\}$ be an orthonormal basis for ${\rm ran}(q)$.  Then if $\|p(a_n)\|<2^{-n}$ for all sufficiently large $n$, $pq$ is compact.
\end{lemma}
\begin{proof}
Let $\epsilon>0$ and let $q_N$ be the projection onto $\overline{{\rm span}}\{a_n\}_{n>N}$ for each $N$.  If $N$ is sufficiently large, then if $v=\sum c_n a_n\in{\rm ran}(q_N)$ is a unit vector,
  \[\|p(v)\|\leq\sum_{n>N}|c_n| \|p(a_n)\|<\sum_{n>N} \frac{1}{2^n}=\frac{1}{2^N}.\]
Thus if we choose $N$ sufficiently large, 
  \[\|\pi(pq)\|=\|\pi(pq_N)\|\leq\|pq_N\|<\epsilon.\]
Since $\epsilon$ is arbitrary, $pq$ is compact.
\end{proof}

\begin{theorem}\label{mathm}
Assume MA.  Let $\Ac\subset\Bc$ be a set of projections of cardinality $<\cc$, and suppose that there is no finite subset $F\subs\Ac$ such that $\bigvee_{p\in F} \pi(p)=1$.  Then there exists a noncompact projection $p\in\Bc$ that is almost orthogonal to every element of $\Ac$.
\end{theorem}
\begin{proof}
Let $Q$ be a countable dense subset of the unit sphere $S\subset\Hc$ with the following property: if $F\subset Q$ is finite, then $Q\cap F^\perp$ is dense in $F^\perp\cap S$.  To obtain such a $Q$, first let $Q_0$ be any countable dense subset of $S$.  For each finite subset $F$ of $Q_0$, enlarge $Q_0$ such that the desired property holds for $F$, and call the result $Q_1$.  Similarly enlarge $Q_1$ to $Q_2$ and so on; we may then set $Q=\bigcup Q_n$.

Now define 
  \[\textnormal{$P=\{(s,F): s$ is a finite sequence of orthogonal elements of $Q$ and $F\subset\Ac$ is finite$\}$}.\]
Order $P$ by saying that $(s',F')\leq(s,F)$ if $s'\supseteq s$, $F'\supseteq F$, and for $n\in{\rm dom}(s')-{\rm dom}(s)$ and $p\in F$, $\|p(s'(n))\|<2^{-n}$.  Then $P$ is clearly ccc since there are only countably many possible values for $s$.  For each $p\in\Ac$, $D_p=\{(s,F)\in P:p\in F\}$ is trivially dense.

We claim that for each $N\in\omega$, $E_N=\{(s,F)\in P:N\subs{\rm dom}(s)\}$ is also dense.  It clearly suffices to show that if $(s,F)\in P$, then there exists $v\in Q$ such that $(s\hat{\;}v,F)<(s,F)$.  Let $n={\rm dom}(s)$ and let $q\in\Bc$ be such that $\pi(q)=1-\bigvee_{p\in F} \pi(p)$.  Let $\{a_n\}$ be an orthonormal basis for ${\rm ran}(q)$; since by hypothesis $q$ is not compact, $\{a_n\}$ is infinite.  Since $q$ is almost orthogonal to each $p\in F$, there is some $m$ such that for each $p \in F$, there are no $2^{-n-1}$-blocks for $p$ above $m$.  Let $q'$ be the projection onto $\overline{{\rm span}}\{a_n\}_{n>m}$ and let $u\in {\rm ran}(q')$ be a unit vector that is orthogonal to $q'(s(k))$ for each $k<n$.  Now let $v\in Q$ be orthogonal to each $s(k)$ such that $\|u-v\|<2^{-n-1}$.  For any $p\in F$, by our choice of $m$, $|\langle u,w\rangle|\leq 2^{-n-1}$ for any unit vector $w\in{\rm ran}(p)$, so
  \[\|p(v)\|\leq \|p(v-u)\|+\|p(u)\|<2^{-n-1}+2^{-n-1}=2^{-n}.\]
Hence $(s\hat{\;}v,F)<(s,F)$.

Since $|\Ac|<\cc$, by MA there is a filter $G\subset P$ meeting each $D_p$ and each $E_N$.  By Lemma \ref{malem}, the projection onto the span of $\bigcup_{(s,F)\in G} s$ is almost orthogonal to each $p\in\Ac$, and is not compact since $G$ meets each $E_N$ so the range is infinite-dimensional.  Thus there is a noncompact projection that is almost orthogonal to every element of $\Ac$, as desired.
\end{proof}

\begin{cor}$\;$
\begin{itemize}
\item[(a):]MA implies all (infinite) maximal aofs have cardinality $\cc$
\item[(b):]MA implies $\kappa$-limits exist in $\Pc$ only for $\kappa=\cc$
\end{itemize}
\end{cor}

\begin{proof}
For (a), if $\{\pi(p_\alpha)\}$ is an infinite aof of cardinality $<\cc$, apply the theorem to $\Ac=\{p_\alpha\}$ to show it is non-maximimal.  For (b), if $((\pi(p_\alpha)),1)$ is a $\lambda$-limit for $\lambda<\cc$, apply the theorem to obtain $p$ that is almost orthogonal to every element of $\Ac=\{p_\alpha\}$; then $\pi(1-p)$ interpolates the limit to give a contradiction.
\end{proof}

\bibliographystyle{amsplain}
.

\end{document}